\documentclass[12pt,english]{article}
\usepackage{amsfonts}
\usepackage{cite,amsfonts,amssymb} 

\def\cD{{\cal D}}
\def\cF{{\cal F}}

\def\cL{{\cal L}}

\def\cR{{\cal R}}

\def\cM{{\cal M}}

\def\sup{\mathop{\rm sup}}

\def\proof{\noindent \medskip {\bf Proof:}$\;\;$}

\newtheorem{prop}{Proposition}

\newtheorem{theor}{Theorem}

\newtheorem{cor}{Corollary}

\newtheorem{rem}{Remark}

\newtheorem{defi}{Definition}

\def\today{}

\title{\bf Stochastic Volterra convolution with L\'evy process}

\author{\large\sf Anna Karczewska \\
 \\ Institute of Mathematics, University of Zielona G\'ora\\
 ul. Podg\'orna 50, 65-246 Zielona G\'ora, Poland\\
 e-mail: A.Karczewska@im.uz.zgora.pl\\}
\date{\today}

\begin{document}

\maketitle

\def\thefootnote{}


\begin{abstract}
In the paper we study 
stochastic convolution appearing in Volterra equation driven by 
so called L\'evy process. By L\'evy process we mean a
process with homogeneous independent increments, 
continuous in probability and cadlag.
\end{abstract}

\section{Introduction}\label{s1}

Let $H$ be a real separable Hilbert space with an inner product
$\langle \cdot ,\cdot\rangle_H$ and a norm $|\cdot|_H$. 
In the paper we consider a stochastic
version of linear, scalar type Volterra equation in $H$ of the form
\begin{equation}\label{E1}
u(t) = \int_0^t \,a(t-\tau)\, Au(\tau)\,d\tau + x + g(t), \quad t\geq 0,
\end{equation}
where $a\in L^1_{\mathrm{loc}}(\mathbb{R}_+)$, $A$ is an 
unbounded linear operator in $H$
with a dense domain $\cD(A)$,
$g$ is an $H$-valued mapping and $x\in H$.

The linear integral equation (\ref{E1}) is a 
subject of many papers connected with applications in different
fields. Among others, the equation (\ref{E1}) 
may be applied to several problems arising in mathematical
physics. For instance,  theory of viscoelasticity provides
numerous problems leading to the Volterra equation of the form (\ref{E1}) 
(see \cite{Pr},  for survey).

We assume that the equation (\ref{E1}) is {\em well-posed} and denote by 
$\cR(t), t\geq 0$, the family of resolvent operators corresponding to (\ref{E1}).
Operators $\cR(t)$ are linear for each $t\geq 0$, 
uniformly bounded on compact intervals,
$\cR(0) x=x$ holds on $\cD(A)$, and $\cR(t)x$ is continuous on 
$\mathbb{R}_+$ for each $x\in \cD(A)$.
Additionally, the following {\em resolvent equation} holds
$$
 \cR(t)\,x = x + \int_0^t a(t-\tau)A\cR(\tau)x\,d\tau
$$
for all $x\in\cD(A)$, $t\geq 0$. 	
For more details concerning well-posedness and resolvent
operators we refer again to the monograph \cite{Pr}.

In the paper we study  equation (\ref{E1}) with an external force $g(t)=Z(t), 
t\geq 0$, where $Z$ is a L\'evy process defined on a stochastic basis 
$(\Omega,\cF ,\cF_t ,P)$. This way the traditional Gaussian framework, when a
Wiener process is the external noise, is extended. Our considerations are
motivated by the growing interest in L\'evy processes in applications, when the
empirical observations simply cannot be explained by means of the Gaussian
distribution. 

So, we arrive at the equation
\begin{equation}\label{E2}
X(t) = \int_0^t \,a(t-\tau)\, AX(\tau)\,d\tau + X_0 + Z(t), 
 \quad\mbox{where~~}  X_0\in H.
\end{equation}

Stochastic equations of Volterra type driven by semimartinagales have
been of course already studied by some authors, for instance \cite{Pro1}
or \cite{Tu1,Tu2}. But our paper treats the subject 
in a different spirit. We use the resolvent operators of the
equation considered, then this way we try to extend the semigroup approach
to the equation (\ref{E2}). For defining the stochastic
convolution we do not use the general semimartingales technique but 
a simpler method.  

We have assumed that $A$ is a closed linear operator in $H$ with the dense 
domain $\cD(A)$. For such a class of operators exists the family 
$\cR(t),\; t\geq 0$, of resolvent operators for (\ref{E2}) 
which fulfills some useful properties (see again \cite{Pr}, Chapter 1). 
Moreover, the family $\cR(t),\; t\geq 0$, is a good enough class of 
operators to be integrands in stochastic integral with respect to 
a process with independent increments. 


Now, we introduce definitions of solutions to (\ref{E2}), analogously
like in previously considered stochastic cases.
\begin{defi} \label{de1}
 An $H$-valued predictable process $X(t)$, $t\in[0,T]$,
 is said to be a {\tt weak solution}
 to (\ref{E2}), if $P(\int_0^t |a(t-\tau)X(\tau)|_H\,d\tau<+\infty)=1$
 and if for all $\xi\in\cD(A^*)$ and all $t\geq 0$ the following 
 equation holds
 $$
  \langle X(t),\xi\rangle_H = \langle X_0,\xi\rangle_H +
  \langle \int_0^t a(t-\tau)X(\tau)d\tau,A^*\xi\rangle_H +
  \langle Z(t),\xi\rangle_H\,, \quad P-a.s.
 $$
\end{defi}
\begin{defi} \label{de2}
  An $H$-valued predictable process $X(t)$, $t\in[0,T]$, 
  is said to be a {\tt mild solution}
 to (\ref{E2}), if $P(\int_0^T |X(\tau)|_H d\tau<+\infty)=1$ and,
 for arbitrary $t\in[0,T]$,
 \begin{equation} \label{E4}
 X(t) = \cR(t)X_0 + \int_0^t \cR(t-\tau)\,dZ(\tau)\,,
 \end{equation}
 where $\cR(t)$ is the resolvent for the equation (\ref{E2}).
\end{defi}

The aim of the paper is to study process called 
{\em stochastic convolution} 
\begin{equation}\label{E5}
Z_{\cR}(t) := \int_0^t \cR(t-\tau)\,dZ(\tau)\,, \quad t\geq 0\;,
\end{equation}
which is the crucial part of the solution (\ref{E4}).

In section 2 we recall the rigorous definition of 
stochastic convolution with L\'evy process.
Then we adapt it
to the Volterra equation (\ref{E2}) driven by L\'evy process
 and use some properties of such integral.
Next, in section 3 we consider particular Volterra equations.

\section{Stochastic convolution}\label{s2}

Assume that $Z(t),\; t\geq 0$, is an $H$-valued process with homogeneous
independent increments (that is, L\'evy process), defined on a fixed probability
space $(\Omega,\cF,\cF_t,P)$, continuous in probability, cadlag and 
$Z_0=0$. It is of great importance in the study of linear and nonlinear
stochastic Volterra equations driven by L\'evy processes to establish first the
basic properties of the process
$ Z_{\cR}(t) = \int_0^t \cR(t-\tau)\,dZ(\tau)$.

The stochastic convolution (\ref{E5}), where $Z$ is a L\'evy 
process, may be defined analogously like stochasic integral in the paper 
\cite{Ch}, that is, as a limit in probability of Stieltjes sums.
This integral coincides with the integral defined by the semimartingales
technique (see \cite{Me} or \cite{Pro2}), but approch used in \cite{Ch}
provides immediately some useful properties of the integral.
The most important is that we obtain 
the explicit formula for the characteristic form 
of the convolution (\ref{E5}).

First we recall some facts concerning stochastic integral with respect to 
L\'evy process used in the paper.

The class $\cL^2_{[u,w]}(H,G)$ of integrands is defined as
follows:\\
$ \Phi :[u,w] \longrightarrow L(H,G), \mbox{~where~} u,w\in \mathbb{R}_+ 
\mbox{~and~}
 H,G \mbox{~are Hilbert spaces},$
such that:\\
~~1) for any $h\in H$, \quad $\Phi h : [u,w] \longrightarrow G$ is measurable,\\
~2) $\int_{u}^w ||\Phi(s)||^2 ds < +\infty$, where $||\cdot||$ means
the operator norm.\\
(In the above definition $L(H,G)$ denote the space of linear bounded 
operators acting from $H$ into $G$.)

\begin{theor}\label{t1} (Theorem 1.3, \cite{Ch})\\
Let $H$ and $G$ be Hilbert spaces and  
a function $\Phi$ belong to the class $\cL^2_{[u,w]}(H,G)$.
Assume that there exists a
sequence $\{ \Phi_n\}$ of step functions that:
\begin{itemize}
\item
for any $n\in N$, $\Phi_n:[u,w]\rightarrow L(H,G)$ and for partition 
$u=s_0<s_1<\ldots <s_n=t$, $ \Phi_n(s)= \Phi_n^k$, where $s\in (s_,s_{k+1}]$,
$k=0,1,\ldots ,n-1;$
\item 
the following condition holds
\begin{equation}\label{E6}
\Phi(s) h = \lim_{n\rightarrow\infty} \Phi_n(s) h \quad \mbox{for any } h\in H
\mbox{ for a.a. } s\in [u,w];
\end{equation}
\item
there exists a function $g\in L^1([u,w])$ that
\begin{equation}\label{E7}
\sup_n ||\Phi_n(s) ||^2 \leq g(s) \quad \mbox{for a.a. } s\in [u,w].
\end{equation}
\end{itemize}
Then the sequence of random variables 
$$ J(\Phi_n) := \int_{(u,w]}  \Phi_n(s)\, dZ_s = \sum_{k=0}^{n-1}
 \Phi_n^k(Z_{s_{k+1}}-Z_{s_{k}})$$
converges in probability and the limit does not depend on the choice of the
sequence $\{ \Phi_n\}$.
\end{theor}
\begin{theor}\label{t2} (Theorem 1.8, \cite{Ch})\\
Let the function $\Phi$ and the sequence $\{\Phi_n\}$ 
satisfy the assumptions of Theorem \ref{t1}.
Then the integral defined as
\begin{equation}\label{E8}
\int_{(u,w]}  \Phi(s)\, dZ_s := P- \lim_{n\rightarrow\infty} 
  \int_{(u,w]}  \Phi_n(s) \,dZ_s 
\end{equation}
is well-defined $G$--valued random variable which has infinitely divisible
distribution.
\end{theor}
\begin{rem}\label{r1}
The integral $\int_u^\infty \Phi(s)\,dZ_s$ is defined as the limit in
probability, as $w\rightarrow\infty$, of the integrals 
$\int_{(u,w]} \Phi(s)\,dZ_s$. The integrals $\int_{-\infty}^t \Phi(s)\,dZ_s$
and $\int_{-\infty}^{+\infty} \Phi(s)\,dZ_s$ are defined analogously.
\end{rem}

As we have already written, the aim of the paper is to study the stochastic
convolution (\ref{E5}), where $Z(t),\; t\geq 0$, is a L\'evy process and 
$\cR (t),\; t\geq 0$, are the resolvent operators to the Volterra equation
(\ref{E2}).
Basing on properties of resolvent operators we can see that the operators 
$\cR (t),\;\, t \in [u,w]$, belong to the class $\cL^2_{[u,w]}(H,H)$.
Actually,  $\cR (t)$ are linear and bounded 
for each $t\geq 0$ and $\cR (t)x$ 
is continuous on $\mathbb{R}_+$ for each $x$ belonging to the 
domain $\cD (A)$ of the operator $A$. Moreover, the function 
$\cR (\cdot)x$, $x\in \cD (A)$, is measurable.
Let $\cR_n (t)$,  $t \in [u,w]$, be step functions defined as follows
\begin{equation}\label{E9}
 \cR_n := \sum_{i=1}^n \cR(s_i)\,\chi_{[t_{i-1},t_i]}\;,
\end{equation}
where $t_0=u,\; t_0<t_1<\ldots <t_{i-1}<t_i<\ldots<t_n=w$,
and $s_i$ is a point from $[t_{i-1},t_i]$.

Let us notice that the 
 functions $\cR_n$  defined by (\ref{E9}) on the interval $[u,w]$,
 satisfy conditions 
  (\ref{E6}) and (\ref{E7}). 
Indeed,  
the interval $[u,w]$ is a compact set in $\mathbb{R}$ and 
the operator $\cR (t)$, $t\in [u,w]$ is continuous with respect to $t$. 
Then the sequence $(\cR_n),\; n\in N$,
of step functions (\ref{E9}) is uniformly convergent to the function $\cR$.
Additionally $(\cR_n),\; n\in N$, are bounded.

So, we may define the stochastic convolution (\ref{E5}) like the stochastic
integral (\ref{E8}), that is, like the limit in probability of integrals of step
functions. Hence, the following theorem comes directly from Theorem \ref{t2}.
\begin{theor}\label{t3}
 Let $\cR (t),\; t\geq 0$, be the family of resolvent operators of the Volterra
 equation (\ref{E2}) and $Z(t),\;\, t\geq 0$, 
 be a L\'evy process. Then the integral 
$$ \int_0^t \cR(t-\tau)\,dZ(\tau) := P - \lim_{n\rightarrow +\infty}
 \int_0^t \cR_n(t-\tau)\,dZ(\tau) $$ 
is well-defined $H$-valued 
random variable which has infinitely divisible distribution. 
\end{theor}

Let us recall that the process $Z(t),\;\, t\geq 0$, as the process 
with independent increments, has the following representation 
(see for instance \cite{Be,GiSk,Ka} or \cite{Sk})
$$ Z_t = at + W_t + \Delta_t \;,$$
where $a\in H, (W_t), t\geq 0$, is an $H$-valued Wiener process and 
$(\Delta_t), t\geq 0 $, is a jump process independent of $(W_t)$. 
This decomposition is clearly unique.
Moreover, for any $t\geq 0 $, the random variable $Z_t$ has L\'evy 
characterization $[ta, t\Theta, tM]$, where $\Theta$ is the covariance
operator of $W_1$ and $M$ is the  L\'evy spectral measure of $\Delta_1$.
This is a consequence of the fact that any infinitely divisible probability 
measure can be viewed as the distribution of a L\'evy process evaluated at 
time 1 and vice versa. Particularly, the famous L\'evy-Khintchine formula
determines the class of characteristic functions corresponding to 
infinitely divisible laws.

Now, we are ready to characterize the convolution (\ref{E5}) as follows.

\begin{theor}\label{t4}
Let L\'evy process $Z(t), \;t\geq 0$, be such that every 
random variable $Z_t$ has L\'evy 
characterization $[ta, t\Theta, tM]$. Then the stochastic convolution\\
$Z_{\cR}(t) = \int_0^t \cR(t-\tau)\,dZ(\tau),~~t\geq 0$,
where $\cR(t)$ are resolvent operators to the Volterra equation (\ref{E2}),
has the following L\'evy  
characterization $[\alpha,Q,\cM]$:  
\begin{eqnarray} \label{E10}
 \alpha & = & \int_0^t \cR(t-\tau)a d\tau + \int_0^t \int_H \cR(t-\tau)x 
[\mbox{\bf 1}_{\{|\cR(t-\tau)\,x|<1\}}-\mbox{\bf 1}_{\{|x|<1\}}] M(dx) d\tau \;;
\nonumber \\[1mm]
  Q & = & \int_0^t \cR(t-\tau)\Theta \cR^*(t-\tau) d\tau \;;\\[1mm]
  \cM  & = &\int_0^t M(\Phi^{-1}(t-\tau)dx)d\tau \;.\nonumber
\end{eqnarray}
\end{theor}

\proof{By Theorem \ref{t3}, the stochastic convolution given by the formula
 $Z_{\cR}(t) = \int_0^t \cR(t-\tau)\,dZ(\tau)\,, \quad t\geq 0$, has 
 infinitely divisible distribution. In order to provide the L\'evy
 characterization of $Z_{\cR}(t)$ it is enough to write the characteristic
 functional of the law of $Z_{\cR}(t)$ and next use the L\'evy-Khintchine 
 formula for the corresponding characteristic exponent. The proof of 
 the theorem is analogous to the proofs of Lemma 1.5 and Theorem 1.8 from
 the paper \cite{Ch}. In our case, the characteristic 
 functional is \\
 $f(y):=\exp \int_0^t \phi(\cR(t-s)y)ds$, for $y\in H$, where 
 $\phi(w)= \log \mathbb{E}(\exp i\langle w,Z_1\rangle)$, $w\in H$. 
{~}  \hfill $\square$}

\begin{cor} \label{C1}
The stochastic convolution (\ref{E5}) is stochastically continuous and
then has a predictable version.
\end{cor}

\begin{theor} \label{t5}
Assume that the operators $\cR(t),~t\geq 0$, are as above. 
Then the stochastic Volterra equation (\ref{E2}) has exactly 
one mild solution.
\end{theor}

\noindent\underline{Comment:} Theorem \ref{t5} comes from uniqueness 
of the resolvent $\cR(t),~t\geq 0$, for deterministic Volterra 
equation (\ref{E1}) and the existence of mild solution to (\ref{E1}).

\begin{theor} \label{t6}
Assume that the operators $A,\cR(t)$ and
the process $Z(t),~t\geq 0$, are like above and the function 
$a\in W_\mathrm{loc}^{1,1}(\mathbb{R}^+)$. 
Let $X$ be an $H$-valued
predictable process with integrable trajectories. If for any $t\in [0,T]$ 
and $\xi\in\cD(A^*)$ the equality 
\begin{equation} \label{E11}
 \langle X(t),\xi\rangle_H = 
 \int_0^t \langle a(t-\tau)X(\tau),A^*\xi\rangle_H\,d\tau +
 \int_0^t \langle \xi,dZ(\tau)\rangle_H
\end{equation}
holds, then 
\begin{equation} \label{E11a}
  X(\cdot)=Z_R(\cdot).
\end{equation}
\end{theor}

\proof{The idea of the proof is the following. First, we prove that if 
(\ref{E11}) is satisfied, then for any 
$\widetilde{\xi}\in C^1([0,T],\cD(A^*) )$ and
$t\in[0,T]$, the following equality holds
\begin{eqnarray} \label{E11b}
 \langle X(t),\widetilde{\xi}(t)\rangle_H & = & 
 \int_0^t \langle (\dot{a}\star X)(\tau)+a(0)X(\tau),
 A^*\widetilde{\xi}\rangle_H\,d\tau +
 \int_0^t \langle \widetilde{\xi}(\tau),dZ(\tau)\rangle_H 
  \nonumber\\  &+&
 \int_0^t \langle X(\tau),\dot{\widetilde{\xi}}(\tau)\rangle_H, 
 \quad\quad P-a.s.
\end{eqnarray}
Next, we take $\xi(\tau):=\cR^*(t-\tau)\xi$ for $\tau\in [0,t]$ and
rewrite (\ref{E11b}).

Then, using properties of resolvent operators, particularly the resolvent
equation, we obtain thesis (\ref{E11a}).
{~}  \hfill $\square$}\\[2mm]

\noindent\underline{Comment:} The above 
Theorem \ref{t6} says that a weak solution to  
(\ref{E2}) is a mild solution to (\ref{E2}).

\begin{theor} \label{t7}
 The stochastic convolution  $Z_R(t)=\int_0^t \cR (t-\tau) dZ(\tau),
 ~t\in [0,T]$, fulfills the equation  (\ref{E11}). 
\end{theor}
\proof{Let us notice that the process $Z_R$ has integrable trajectories
(the set of ,,discontinuity" is at most countable).\\
For any $\xi\in\cD(A^*)$ we may write 
\begin{eqnarray*} 
\int_0^t \langle a(t-\tau) Z_R(\tau),A^*\xi \rangle_H d\tau 
\!&\! = \!&\! ~\\
  \mbox{(from~(\ref{E5}))} \!&\! = \!&\!
\int_0^t \langle a(t-\tau)\int_0^\tau\cR(\tau-\sigma)dZ(\sigma),A^*\xi 
 \rangle_H \\
  \mbox{(from~Dirichlet~formula}  \!&\! \mbox{and} \!&\!
  \mbox{stochastic~Fubini~theorem)} \\ \!&\! = \!&\!
\int_0^t \langle \left[ 
\int_0^t a(t-\tau)\cR(\tau-\sigma)d\tau\right]
 dZ(\sigma),A^*\xi \rangle_H \\
  \!&\! = \!&\!   
  \langle \int_0^t \left[ \int_0^{t-\sigma} \!\!\! a(t-\sigma-z)
  \cR(z)dz\right]dZ(\sigma) ,A^*\xi \rangle_H  \\ \!&\! = \!&\!
  \langle \int_0^t A[(a\star \cR)(t-\sigma)]dZ(\sigma),\xi \rangle_H \\ 
\mbox{(from~resolvent~equation)}  \!&\! = \!&\!
 \langle \int_0^t [\cR(t-\sigma)-I]dZ(\sigma),\xi \rangle_H \\\!&\! = \!&\!
 \langle \int_0^t \cR(t-\sigma)dZ(\sigma),\xi \rangle_H -
 \langle \int_0^t  dZ(\sigma),\xi \rangle_H \,.
\end{eqnarray*}
Hence, we obtained the following equation
$$
 \langle Z_R(t),\xi \rangle_H \ = \int_0^t \langle a(t-\tau) Z_R(\tau),
 A^*\xi \rangle_H d\tau + \int_0^t \langle \xi,dZ(\tau)\rangle_H
$$ 
for any $\xi\in\cD(A^*)$.
{~}  \hfill $\square$}

\begin{cor} \label{c2}
 Assume that the operator $A$ is bounded. Then 
 $$
 Z_R(t) = \int_0^t  a(t-\tau) AZ_R(\tau)d\tau
  + \int_0^t dZ(\tau) \,.
$$ 
\end{cor}

\section{Particular cases}

Let us notice that till now we have not assumed that the resolvent operators 
$\cR (t),\; t\geq 0$, to the Volterra equation (\ref{E2}) have bounded
variation. The stochastic integral with respect to L\'evy process and the
stochastic convolution (\ref{E5}) have been defined and characterized 
(formula (\ref{E10})) without this 
assumption. Hence, if the stochastic Volterra equation (\ref{E2}) driven by 
L\'evy process $Z(t),\; t\geq 0$, is well-posed, that is,
the equation (\ref{E2})  has the resolvents
$\cR (t),\; t\geq 0$, then Theorem \ref{t4} for the convolution  
$Z_R(t),\; t\geq 0$, holds. Besides, there are some Volterra equations 
which admit 
resolvent operators with bounded variation. These equations may be treated in
a different way. So, in this section we study such particular case of 
Volterra equations.

We will take the following common assumption.

\noindent {\bf Assumption} (A) 
There exists the family 
$\cR(t),\; t\geq 0$, of resolvent operators to the 
Volterra equation (\ref{E2}), the operators have bounded variation and 
$Z(t),\; t\geq 0$, is an $H$-valued L\'evy process.

If the resolvent operators  $\cR(t),\; t\geq 0$, have bounded variation, we may
use for stochastic convolution (\ref{E5}) the classical integration by parts.
Then we may write
\begin{equation} \label{E12}
\int_{(a,b]} \cR(t-\tau) dZ(\tau) = \cR(t-b) Z(b) - \cR(t-a) Z(a)
  - \int_{(a,b]}  Z(\tau-) d\cR(t-\tau) \;,
\end{equation} 
where $Z(t),\; t\geq 0$, is a stochastic  L\'evy process.

The properties of the stochastic integral $\int_{(a,b]} \cR(t-\tau) dZ(\tau)$
may be obtained from the right hand side of (\ref{E12}).
Among others, we can deduce the following results.

\begin{prop}\label{p1}
Under the assumption (A) we have 
$$ \langle x^*,\int_{(a,t]} \cR(t-\tau) dZ(\tau) \rangle = 
 \int_{(a,t]} \cR(t-\tau) d\langle x^*,Z(\tau)\rangle $$
 for $ x^* \in H^*$, where $ H^*$ is the dual space to $H$.
\end{prop}
\vspace{2mm}

\begin{prop}\label{p2}
If the assumption (A) holds, the function 
$$ t\longrightarrow \int_{(a,t]} \cR(t-\tau) dZ(\tau) $$
 is $H$-valued random variable with infinite divisible distribution.
\end{prop}

This fact follows from (\ref{E12}) and the approximation by 
Riemann-Stieltjes sums. See, e.g. \cite{JuVe}.

\begin{prop}\label{p3}
Let the assumption (A) be satisfied. Then 
\begin{equation}\label{E13}
 \log\left[ \hat{\cL} (\int_{(0,t]} \cR(t-\tau) dZ(\tau))(\lambda)\right]
=  \int_{(0,t]} \log [\hat{\cL}(Z(1))(\lambda\cR(t-s)) ] ds\;,
\end{equation}
 for $\lambda \in H^*$, where $H^*$ is the dual space to $H$ and $\hat{\cL}$ 
 denotes the characteristic functional of the probability distribution of
 the appropriate random variable.
\end{prop}

\proof{The above formula (\ref{E13}) comes from Lemma 1.1, \cite{JuVe}
and the definition of the convolution (\ref{E5}).
\hfill $\square$
}

\begin{prop}\label{p4}
If the assumption (A) holds, then
$ \int_{(0,t]} \cR(t-\tau) dZ(\tau)$ converges in norm w.p. 1, as 
$t\rightarrow +\infty$, if and only if $ \int_{(0,t]} \cR(t-\tau) dZ(\tau)$ 
converges in the distribution as $t\rightarrow \infty$.
\end{prop}

\proof{Proposition comes from Lemma 1.2, \cite{JuVe}. 
\hfill $\square$
}
\vspace{2mm}

We finish the paper by the examples of Volterra equations. 

Let us consider the case when the function $a$ is
completely positive. This class of  kernels is very important
in the theory of Volterra equations and arises naturally in applications. 
If the function $a\in L^1_{\mathrm{loc}}(\mathbb{R}_+) $ is completely positive,
then $s(\cdot ,\gamma)$, the solution to the integral equation
\begin{equation}\label{E3}
s(t) +\gamma \int_0^t \,a(t-\tau)\, s(\tau)\,d\tau = 1, \quad t\geq 0\;,
\end{equation}
is nonnegative and nonincreasing for any $\gamma > 0$.
More precisely, under this condition $s(t)\in [0,1]$.
There is a relationship between the resolvent operator $\cR(t)$ to the equation
(\ref{E2}) and the corresponding function $s(t,\gamma)$ fulfilling 
(\ref{E3}). Namely, if $-\gamma$ is an eigenvalue of $A$ with eigenvector
$z\neq 0$, then $\cR(t)z = s(t,\gamma)z,\: t\geq 0$. 

Because every function monotonic on the interval $[0,T]$ has bounded 
variation on $[0,T]$, it is enough to choose any function $a(t),~t\geq 0$,
which is completely positive.

Particularly, 
let $H=L^2(0,1)$ and $Au=D^2 u$ with $\cD(A) = H^2(0,1) \cap H^1_0(0,1)$.
The functions $e_k(\xi)= \sqrt{2/\pi}\sin k\xi,\; \xi\in[0,1],\; k\in N$,
form an orthonormal sequence of eigenfunctions of the operator $A$, 
corresponding to the eigenvalues $-\mu_k = -\pi^2k^2,\; k\in N$. 
When we set additionally, that the function $a(t)=e^{-t}~$ for $t\geq 0$,
then we obtain in this case
$$  s(t,\mu)=(1+\mu)^{-1} [1+\mu\,e^{-(1+\mu)t}],\quad t,\mu>0\;. $$

There exists the resolvent $\cR(t),\; t\geq 0$, to the equation (\ref{E2}) in 
this case and is determined by 
$ \cR(t)\,e_k = s(t,\mu_k) e_k, \quad k\in N\;.  $
So, the function $\cR(t),\; t\geq 0$, is monotonic in $[0,1]$, because the
function $s(t,\mu)$ is. Hence, $\cR(t)$ has bounded variation on $[0,1]$.

Other examples of Volterra equations with resolvent operators having bounded
variation may be found in the monograph \cite{Pr}.

\end{document}